\documentclass[12pt]{amsart}
\usepackage[noadjust]{cite}
\usepackage[final]{graphicx}
\usepackage[margin=1in,letterpaper]{geometry}
\usepackage{booktabs}
\usepackage{url}
\usepackage[T1]{fontenc}
\usepackage{lmodern}
\usepackage[utf8]{inputenc}
\usepackage{mathtools}
\usepackage[bookmarks=true,%
    colorlinks=true,%
    linkcolor=blue,%
    citecolor=blue,%
    filecolor=blue,%
    menucolor=blue,%
    urlcolor=blue]{hyperref}
\hypersetup{pdftitle={Characterizing the universal rigidity of generic
    frameworks}}
\hypersetup{pdfauthor={Steven J. Gortler and Dylan P. Thurston}}

\usepackage{hyphenat}

\newread\testin

\def\mathcenter#1{%
  \vcenter{\hbox{$#1$}}%
}

\def\mfigb#1{
        \mathcenter{\includegraphics[trim=-1 -1 -1 -1]{#1}}
}

\graphicspath{{draws/}}

\newtheoremstyle{citing}
  {}
  {}
  {\itshape}
  {}
  {\bfseries}
  {.}
  {.5em}
  {\thmnote{#1 #3}}

\theoremstyle{citing}
\newtheorem*{citingthm}{Theorem}

\theoremstyle{plain}

\newcommand{\RR}{\mathbb R}

\newcommand{\QQ}{\mathbb Q}

\newcommand{\EE}{\mathbb E}

\newcommand{\kk}{\mathbf{k}}

\newcommand{\abs}[1]{{\lvert #1 \rvert}}

\newcommand{\bdy}{\partial}

\DeclareMathOperator{\Eucl}{Eucl}

\DeclareMathOperator{\Aff}{Aff}

\DeclareMathOperator{\tr}{tr}

\DeclareMathOperator{\rank}{rank}

\DeclareMathOperator{\Int}{Int} 
\DeclareMathOperator{\ext}{ext} 

\DeclareMathOperator{\GL}{\textit{GL}}

\theoremstyle{plain}
\numberwithin{equation}{section}
\newtheorem{theorem}{Theorem}
\newtheorem{proposition}[equation]{Proposition}
\newtheorem{otheorem}[equation]{Theorem}
\newtheorem{lemma}[equation]{Lemma}
\newtheorem{corollary}[equation]{Corollary}

\theoremstyle{definition}
\newtheorem{definition}[equation]{Definition}

\newtheorem{question}[equation]{Question}

\theoremstyle{remark}

\newtheorem{remark}[equation]{Remark}

\newcommand{\Edges}{{\mathcal E}}
\newcommand{\Verts}{{\mathcal V}}

\newcommand{\prho}{p}

\begin{document}
\title{Characterizing the universal rigidity of generic frameworks}

\author[Gortler]{Steven J. Gortler}
\address{School of Engineering and Applied Sciences\\
  Harvard University\\
  Cambridge, MA 02138}
\email{sjg@cs.harvard.edu}

\author[Thurston]{Dylan~P. Thurston}
\address{Department of Mathematics\\
         Indiana University\\
         Bloomington IN 47405}
\email{dpthurst@indiana.edu}

\begin{abstract}
  A \emph{framework} is a graph and a map from its vertices to~$\EE^d$
  (for some $d$).  A framework is \emph{universally rigid} if
  any framework in any dimension with the same graph and edge lengths
  is a Euclidean image of it.  We show that a 
generic universally rigid
  framework has a positive semi-definite stress matrix
of  maximal rank.  Connelly showed that the existence of such a positive
  semi-definite stress matrix is sufficient for universal rigidity, so
  this provides a characterization of universal rigidity for generic
  frameworks. We also extend our argument to give a new result on 
the genericity of strict complementarity in semidefinite programming.
\end{abstract}

\date{\today}


\maketitle

\section{Introduction}
\label{sec:intro}

In this paper we characterize generic frameworks which are
universally rigid in $d$\hyp dimensional Euclidean space. 
A framework is universally rigid if, modulo
Euclidean transforms, 
 there is no other framework of the same graph in 
\emph{any} dimension that has the same edge lengths.
A series of papers by Connelly~\cite{Connelly82:RigidityEnergy,
Connelly99:TensegrityStable,
Connelly01:StressStability} and later
Alfakih~\cite{Alfakih07:DimensionalRigidity,alfakih2007universal}
described a 
sufficient condition for a generic 
framework to be universally 
rigid. In this paper we show that this is also a necessary
condition.

Universally rigid frameworks are especially relevant when applying
semidefinite programming techniques to graph embedding problems. 
Suppose some vertices are embedded in $\EE^d$ and we are told the 
distances between some of the pairs of vertices. The graph embedding
problem is to compute the embedding 
(up to an unknown Euclidean transformation)
from the data. This problem is computationally difficult as the
graph embeddablity question is 
in general
NP-HARD~\cite{Saxe79:EmbedGraphsNP}, but
because of its utility, many heuristics have been attempted. 
One approach is to use semidefinite programming to find an 
embedding~\cite{linial1995gga}, 
but such methods are not able to find the solutions that are
specifically $d$-dimensional, rather than embedded in some larger dimensional
space.
However, when the underlying
framework is, in fact,
universally rigid, then the distance data itself automatically
constrains the dimension, and therefore
the correct answer will be (approximately)
found by the semidefinite program.
The connection between rigidity and semi-definite programming
was first explored by So and Ye~\cite{SY07:SemidefProgramming}.

Many rigidity questions are easier to answer in the 
generic case (see Definition~\ref{def:generic} below), where various
singular cases can simply be ignored. It is also a reasonable assumption to
make in many graph embedding problems that arise in practice.
For example, in the problem of 
sensor network localization~\cite{eren2004rigidity}, we assume that there
are a set of robotic sensors moving independently in space, and that
some inter-sensor distances are measured physically. For such a setting,
we  only expect to observe generic behavior.

We also discuss the more general topic of strict complementarity 
in semidefinite programming. Strict complementarity 
is a strong form of duality
that is needed for the fast convergence of many interior point 
optimization algorithms. Using our arguments from universal rigidity,
we show that if the semidefinite program has a sufficiently 
generic primal solution, then it must satisfy strict complementarity. 
This is in distinction with previous results on strict 
complementarity~\cite{alizadeh1997complementarity,pataki2001generic} 
that require the actual parameters of the program to be generic.
In particular, our result applies to programs where the solution is
of a lower rank than would be found generically.

\subsection{Rigidity definitions}
\label{sec:definitions}
\begin{definition}\label{def:config-space}
  A \emph{graph}~$\Gamma$ is a set of $v$ vertices $\Verts(\Gamma)$
  and $e$ edges~$\Edges(\Gamma)$, where $\Edges(\Gamma)$ is a set of\
  two\hyp element subsets of $\Verts(\Gamma)$.  We will typically drop
  the graph~$\Gamma$ from this notation.
  A \emph{configuration}~$\prho$  is 
  a mapping from $\Verts$ to~$\EE^v$.
Let $C(\Verts)$
be the space of configurations.
For $\prho\in C(\Verts)$ and $u
  \in \Verts$, 
let $\prho(u)$ denote the image of $u$ under~$\prho$.  
  Let $C^d(\Verts)$
denote the space of configurations that lie entirely 
in the $\EE^d$, contained as the first $d$ dimensions
of~$\EE^v$.
A \emph{framework} $(\prho,\Gamma)$ is the pair of a graph and a configuration
of its vertices.
For a given graph~$\Gamma$
  the \emph{length-squared function}
  $\ell_\Gamma:C(\Verts)\rightarrow\RR^e$ is the function assigning to each
  edge of $\Gamma$ its squared edge length in the framework.  That is,
  the component of $\ell_\Gamma(\prho)$ in the direction of an
  edge $\{u,w\}$ is $\abs{\prho(u)-\prho(w)}^2$.
\end{definition}

\begin{definition}
  \label{def:generic}
  A configuration in $C^d(\Verts)$ is \emph{proper} if it
  does not lie in any affine subspace of $\EE^d$ of dimension less
  than~$d$.  It is \emph{generic}
if its first $d$ coordinates (i.e., the coordinates not constrained to
be~$0$) do
  not satisfy any algebraic equation with rational 
  coefficients.
\end{definition}

\begin{remark}
  A generic configuration in $C^d(\Verts)$ with at least $d+1$ vertices is proper.
\end{remark}

\begin{definition}
  \label{def:universally-rigid}
The configurations  $\prho, q$ in $C(\Verts)$ are \emph{congruent}
if they are related by an element of the group of $\Eucl(v)$ of rigid
motions of~$\EE^v$.

A framework $(\prho,\Gamma)$ with $\prho \in C^d(\Verts)$
is \emph{universally rigid} if
  any other configuration in $C(\Verts)$ with 
the same edge lengths under $\ell_{\Gamma}$
is a 
configuration
congruent to $\prho$.

A graph $\Gamma$ is \emph{generically universally rigid} in $\EE^d$
if any generic 
framework $(\prho,\Gamma)$ with $\prho \in C^d(\Verts)$
is universally rigid.
\end{definition}

Essentially, 
universal rigidity means that the lengths of the edges of $\prho$ are 
consistent with essentially only one embedding of $\Gamma$ in
\emph{any} dimension, up to~$v$.  (In general, a configuration in any
higher dimension
can be related by a rigid motion to a configuration in $\EE^v$.)
This is stronger than \emph{global rigidity}, where the lengths fully
determine the embedding in the smaller space
$\EE^d$. And global rigidity is, in turn,
stronger than (local) \emph{rigidity}, which only rules out continuous
flexes in $\EE^d$
that preserve edge lengths.

\begin{remark}
  Universal rigidity or closely related notions have also been called
  \emph{dimensional rigidity}~\cite{Alfakih07:DimensionalRigidity},
  \emph{uniquely localizable}~\cite{SY07:SemidefProgramming}, and
  \emph{super stability}~\cite{Connelly99:TensegrityStable}.
\end{remark}

\begin{remark}
  In Definition~\ref{def:universally-rigid}, it would be equivalent to
  require that $(\prho,\Gamma)$ be (locally) rigid in $\EE^v$ 
  since by a result of Bezdek and
  Connelly~\cite{BC04:KneserPoulsen} any two frameworks with the same
  edge lengths can be connected by a smooth path in a sufficiently
  large dimension.
\end{remark}

\begin{remark}
  Universal rigidity, unlike local and global rigidity, is not a generic
  property of a graph: for many graphs, some generic frameworks are universally
  rigid and others are not.  For instance, an embedding of a 4-cycle
  in the line~$\EE^1$ is universally rigid iff one side is long in the
  sense that its length
  is equal to the sum of the lengths of the others.
  On the other hand, some graphs, such as a simplex,
and any trilateration  graph~\cite{eren2004rigidity}
are generically universally rigid in $\EE^d$. 
(A $d$-trilateration graph has the property that one can order its vertices
such that the following holds:
the first $d+2$ vertices are part of a simplex in $\Gamma$ and
each subsequent vertex is adjacent in  $\Gamma$ 
to $d+1$ previous vertices.)
We do not know a
 characterization of  generically universally rigid graphs.
\end{remark}

\begin{figure}[t]
\includegraphics{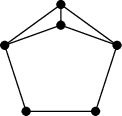}
\qquad
\includegraphics{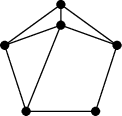}
\qquad
\includegraphics{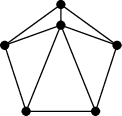}
\qquad
\includegraphics{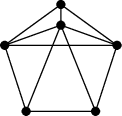}
\qquad
\includegraphics{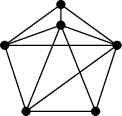}
\caption{Planar graph embeddings with
increasing types of rigidity.  From left to right, 
generically locally flexible graph,
generically locally rigid graph,
generically globally rigid graph, 
generic universally rigid framework, 
generically universally rigid graph.
(More precisely, the fourth framework is not itself generic, but it
and any generic
framework close to it are
universally rigid.)
This paper focuses on frameworks of the type of the fourth one.
The fifth graph is generically universally rigid as it is a 
2-trilateration graph.
}
\label{seq}
\end{figure}

See Figure~\ref{seq} for examples of embedded graphs with increasing
types of rigidity.

\subsection{Equilibrium stresses}

\begin{definition}
  \label{def:stress-matrix}
  An  \emph{equilibrium stress matrix} of a framework~$(\prho,\Gamma)$
is a matrix~$\Omega$
  indexed by $\Verts\times \Verts$ so that
  \begin{enumerate}
  \item for all $u,w \in \Verts$, $\Omega(u,w) = \Omega(w,u)$;
  \item for all $u,w \in \Verts$ with $u \ne w$ and $\{u,w\} \not\in \Edges$,
    $\Omega(u,w) = 0$;
  \item \label{item:stress-row-sum} for all $u \in \Verts$, $\sum_{w\in\Verts} \Omega(u,w) = 0$; and
  \item \label{item:stress-equilib} for all
    $u \in \Verts$,
    $\sum_{w\in\Verts} \Omega(u,w)\prho(w) = 0$.
  \end{enumerate}
(The last condition is the equilibrium condition.)
Let $S_\Gamma(\prho)$ be the linear space of equilibrium stress matrices
for $(\prho,\Gamma)$.
\end{definition}

Conditions \eqref{item:stress-row-sum} and \eqref{item:stress-equilib} give us:
\begin{equation}
\forall u \in \Verts \;\;
  \sum_{\{w \in \Verts \mid \{u,w\} \in \Edges\}}\!\!\Omega(u,w)(\prho(w) -
  \prho(u)) = 0.
\end{equation}

The kernel of an equilibrium stress matrix $\Omega$ of a
framework $(\prho,\Gamma)$
always contains the subspace of $\RR^v$
spanned by the coordinates of~$\prho$ along each axis and the
vector~$\vec 1$ of all $1$'s. This corresponds to the fact that any
affine image of $\prho$ satisfies all of the equilibrium stresses in  
$S_\Gamma(\prho)$.
If $\prho$ is a proper $d$-dimensional configuration, these kernel vectors 
span a
$(d+1)$-dimensional space, so for such frameworks
$\rank \Omega \leq v-d-1$.

\begin{definition}
  \label{def:conic}
We say that the edge directions of $(\prho,\Gamma)$ with $\prho \in C^d(\Verts)$
are \emph{on a conic at
  infinity} if there exists a symmetric $d$-by-$d$ matrix~$Q$ such that
for all edges $(u,v)$ of 
$\Gamma$, we have
\[ [\prho(u)-\prho(v)]^t Q [\prho(u)-\prho(v)] = 0,\]
where the square brackets mean the projection of a vector in $\EE^v$
to $\EE^d$ (i.e., dropping the $0$'s at the end of $\prho(u)$).
\end{definition}

\begin{remark}
  If the edges of~$(\prho,\Gamma)$ are not on a conic at infinity in
  $C^d(\Verts)$, then in
  particular $\prho$ is proper.
\end{remark}

A framework
has  edges on a conic at infinity iff  there is
a continuous family of $d$-dimensional non-orthonormal affine transforms that
preserve all of the edge lengths.
However, if the configuration is not proper, this family of affine
transforms might all be congruent to the original one.

\subsection{Results}

Connelly in a series of papers has studied the relationship between
various forms of rigidity and equilibrium stress matrices. In particular he 
proved the
following theorem that gives a sufficient condition for a 
framework to be universally rigid. 

\begin{otheorem}[Connelly]
\label{thm:suff}
Suppose that $\Gamma$ is a graph with $d+2$ or more vertices
and $\prho$ is a (generic or not)
configuration in $C^d(\Verts)$.
Suppose that there is an equilibrium
 stress matrix $\Omega \in S_\Gamma(\prho)$ 
that is positive semidefinite (PSD) and that $\rank \Omega = v-d-1$.
Also suppose that 
the edge directions of $\prho$ do not lie on a conic at infinity. 
Then
$(\prho,\Gamma)$ is universally rigid.
\end{otheorem}

Connelly also proved
a lemma that allows him to ignore the conic at infinity issue
when $\prho$ is generic.

\begin{lemma}[Connelly]
\label{lem:cinf}
Suppose that $\Gamma$ is a graph with $d+2$ or more vertices
and $\prho$ is a generic 
configuration in~$C^d(\Verts)$.
Suppose that there is an equilibrium
 stress matrix $\Omega \in S_\Gamma(\prho)$
with  rank  $v-d-1$.
Then the edge directions of $\prho$ do not lie on a conic
at infinity.
\end{lemma}

Putting these together, we can summarize this as 

\begin{corollary}
\label{cor:suff2}
Suppose that $\Gamma$ is a graph with $d+2$ or more vertices
and $\prho$ is a generic
configuration in $C^d(\Verts)$.
Suppose that there is a PSD equilibrium
 stress matrix $\Omega \in S_\Gamma(\prho)$ 
with $\rank \Omega = v-d-1$.
Then
$(\prho,\Gamma)$ is universally rigid.
\end{corollary}

The basic ideas for the proof of Theorem~\ref{thm:suff} appear 
in~\cite{Connelly82:RigidityEnergy} where they  were applied to show
the universal rigidity of
Cauchy polygons.
It is also described in~\cite{Connelly99:TensegrityStable}
and stated precisely
in~\cite[Theorem~2.6]{Connelly01:StressStability}.
Lemma~\ref{lem:cinf} can be derived from the proof of Theorem 1.3 
in~\cite{Connelly05:GenericGlobalRigidity}.
Corollary~\ref{cor:suff2} in two dimensions
is summarized in 
Jord\'an and Szabadka~\cite{jordan2009operations}.

A different set of
sufficient conditions for the related concept of
dimensional rigidity  were described by
Alfakih~\cite{Alfakih07:DimensionalRigidity}. 
In~\cite{alfakih2007universal} he showed that these conditions
were equivalent to the existence of a PSD equilibrium stress matrix
of maximal rank.
In~\cite{alfakih2007universal} he also 
showed  that a generic dimensionally rigid framework  (with $v \geq d+2$)
must be universally rigid;
this  has the same effect as Lemma~\ref{lem:cinf}, and thus results in 
Corollary~\ref{cor:suff2}. In these papers, he also conjectured
that for a generic universally rigid framework, a maximal rank
PSD
equilibrium stress matrix must exist.

In the related context of frameworks with pinned anchor vertices,
So and Ye~\cite{SY07:SemidefProgramming}
showed the appropriate analogue of Theorem~\ref{thm:suff} follows
from complementarity in semidefinite programming (see 
section~\ref{sec:sdp} below for more on this).

In this paper our main result is the converse to
Corollary~\ref{cor:suff2}:

\begin{theorem}\label{thm:generic-ur-min-stress-kernel}
  A universally rigid framework $(\prho,\Gamma)$, with $\prho$ generic
in $C^d(\Verts)$ and having $d+2$ or more vertices,
has a PSD
equilibrium stress matrix with rank $v-d-1$.
\end{theorem}

\begin{remark}
  Alfakih has given an example \cite[Example
  3.1]{Alfakih07:DimensionalRigidity} showing that
  Theorem~\ref{thm:generic-ur-min-stress-kernel} is false if we drop
  the assumption that $\prho$ is generic.
  For any universally rigid framework (generic or not), it's not hard
  to see that there is a non-zero PSD
  equilibrium stress matrix.  (See \cite[Theorem
  5.1]{Alfakih09:BarSDP}.)  The difficulty in
  Theorem~\ref{thm:generic-ur-min-stress-kernel} is finding a stress
  matrix of high rank.
\end{remark}

Theorems \ref{thm:suff} and \ref{thm:generic-ur-min-stress-kernel}
compare nicely with the situation for global rigidity, where a
framework is generically globally rigid iff it has an equilibrium
stress matrix of rank $v-d-1$ (with no PSD constraint).  Sufficiency
was proved by Connelly\cite{Connelly05:GenericGlobalRigidity}, and
necessity was proved by the authors together with Alex Healy~\cite{GHT10:GGR}.

\begin{question}
  For a graph~$\Gamma$ that is
generic globally rigid $\EE^d$, is there always a generic framework 
in $C^d(\Verts)$
that is
  universally rigid?
\end{question}

In this paper, we also prove more general versions of
Theorem~\ref{thm:generic-ur-min-stress-kernel} in the general context
of convex optimization and
strict complementarity in semidefinite programming.  See
Theorem~\ref{thm:urexp} in Section~\ref{sec:convex} and
Theorem~\ref{thm:sdp-opt} in Section~\ref{sec:sdp}.


\section{Algebraic geometry preliminaries}
\label{sec:algebr-geom}

We start with some preliminaries about semi-algebraic sets from real
algebraic geometry, somewhat specialized to our particular case.  For
a general reference, see, for instance, the
book by Bochnak, Coste, and Roy~\cite{BCR98:RealAlgGeom}.
\begin{definition}
  An affine real \emph{algebraic set} or \emph{variety}~$V$ 
  contained in~$\RR^n$ is a subset of $\RR^n$ that is
  defined by a set of algebraic equations.
  It is \emph{defined over~$\QQ$} if the equations can be taken to
  have rational coefficients.
  An algebraic set has a \emph{dimension}
  $\dim(V)$, which we will define as the largest $t$ for which there
  is an open subset of~$V$ (in the Euclidean topology) homeomorphic to
  $\RR^t$.
\end{definition}

\begin{definition}
  A \emph{semi-algebraic set}~$S$
is a subset of $\RR^n$ defined by algebraic
  equalities and inequalities;
alternatively (by a non-trivial
  theorem), it is the image of an
  algebraic set (defined only by equalities) under an algebraic map.
  It is \emph{defined over~$\QQ$} if the equalities and inequalities
  have rational coefficients.
  Like algebraic sets, a semi-algebraic set has a well defined
  (maximal) dimension~$t$.
\end{definition}

We next define genericity in larger generality and give some basic
properties.

\begin{definition}
  A point in
  a (semi-)algebraic set~$V$ defined over~$\QQ$ is
  \emph{generic} if its
  coordinates do not satisfy
  any algebraic equation with coefficients in~$\QQ$
  besides those that are satisfied by every
  point on~$V$.

  A point $x$ on a 
semi-algebraic set $S$ is \emph{locally generic} if for small enough
$\epsilon$, 
$x$ is generic in $S \cap B_\epsilon(x)$.
\end{definition}
\begin{remark}
A semi-algebraic set $S$ will have no generic points
if its Zariski closure over $\QQ$ is reducible over $\QQ$. In this case
all points in each irreducible component will satisfy some specific equations
not satisfied everywhere over $S$.
But even in this case,
almost every point in the set  will be generic within its 
own  component. Such points will still be locally generic in $S$.
\end{remark}

We will need a few elementary lemmas on generic points in
semi-algebraic sets.
\begin{lemma}
  \label{lem:dense}
  If $X$ and $Y$ are both semi-algebraic sets defined over $\QQ$, 
with
  $X\subset Y$ and $X$ dense in $Y$ (in the Euclidean topology), 
then $X$ contains all of the locally generic points of~$Y$.  
\end{lemma}
\begin{proof}
Due to the density assumption, $Y\backslash X$ must be a semi-algebraic set
defined over $\QQ$
with dimension less than that of $Y$. The Zariski closure of 
a semi-algebraic
set maintains its dimension, thus all points in $Y\backslash X$ must 
satisfy some algebraic equation that is non-zero over $Y$ and thus these
points must be
non-generic.  To see that a point $y \in Y \setminus X$ is not locally
generic either,
apply this argument in an
$\epsilon$-neighborhood of~$y$.  (That is, apply the argument to $Y
\cap B_\epsilon(y)$ and $X \cap B_\epsilon(y)$.)
\end{proof}

\begin{lemma}
  \label{lem:image-generic}
  Let $V$ be a  semi-algebraic set, $f$ be 
an algebraic map from $V$ to a range space~$X$, and $W$\!
be an  semi-algebraic set contained in the image of $f$,
with $V$, $W$, and $f$ all defined over~$\QQ$.
If $x_0\in V$\! is generic and $f(x_0) \in W$\!,
then 
$f(x_0)$ is generic in~$W$\!.
\end{lemma}

\begin{proof}
  Let $\phi$ be any algebraic function on~$X$ with rational coefficients
  so that $\phi(f(x_0)) =
  0$.  Then $\phi \circ f$ vanishes at~$x_0$, which is generic in~$V$, so
  $\phi\circ f$ vanishes identically on~$V\!$.  This implies that $\phi$
  vanishes on~$W$\!.
  Thus any algebraic function defined 
over~$\QQ$ that vanishes at $f(x_0)$ must vanish on all of $W\!$.
This proves that 
$f(x_0)$ is generic in~$W\!$.
\end{proof}


\section{The geometry of PSD stresses}

We now turn to the main construction in our proof.

\subsection{The measurement set}
\begin{definition}
  The $d$-dimensional \emph{measurement set}~$M_d(\Gamma)$ 
  of a graph~$\Gamma$ is defined to be the image of
  $C^d(\Verts)$ under the map $\ell_\Gamma$.  These are nested by
  $M_d(\Gamma) \subset M_{d+1}(\Gamma)$ and eventually stabilize at
  $M_{v-1}(\Gamma)$, also called the
  \emph{absolute measurement set} $M(\Gamma)$.
\end{definition}

Since $M_d(\Gamma)$ is the image of 
$C^d(\Verts)$ under an algebraic map, 
by Lemma~\ref{lem:image-generic}, if $\prho$ is
generic in $C^d(\Verts)$ then $\ell_\Gamma(\prho)$ is generic in $M_d(\Gamma)$
(which must be irreducible).

\begin{lemma}
  The set $M(\Gamma)$ is convex.
\end{lemma}

\begin{proof}
  The squared edge lengths of a framework~$\prho$ are  computed by
  summing the squared edge
  lengths of each coordinate projection of~$\prho$, so $M_d(\Gamma)$ is
  the $d$-fold Minkowski sum of $M_1(\Gamma)$ with itself.  Since
  $M_1(\Gamma)$ is invariant under scaling by positive reals,
  $M_d(\Gamma)$ can also be described as an iterated chord variety,
  and in particular $M(\Gamma)$ is the convex hull.
\end{proof}

A particular case of interest is $M(\Delta_v)$,
the absolute measurement set of the complete graph on the
vertices.
It is a cone in 
$\RR^{\binom{v}{2}}$.

\begin{lemma}\label{lem:meas-sdp}
  The measurement set $M_d(\Delta_v)$ is isomorphic (as a subset of
  the linear space $\RR^{\binom{v}{2}}$) to the
  set of PSD $(v-1)\times (v-1)$ matrices of rank at most~$d$.
\end{lemma}

\begin{proof}
There is a standard linear map from $M(\Delta_v)$ to the convex cone of
$v-1$ by $v-1$ positive semidefinite matrices, which we recall for
the reader.
(See,
e.g., \cite{schoenberg1935remarks,gower1985properties}.)
For a point $x \in M_d(\Delta_v)$, 
by the universal rigidity of the simplex, up to $d$-dimensional
Euclidean isometry
there is a
unique~$\prho\in
C^d(\Verts)$
with $x= \ell_\Delta(\prho)$.
If we additionally constrain the last vertex
to be at the origin, then $\prho$ is unique up to
an orthonormal transform.
Think of such a constrained~$\prho$ as a
$(v-1) \times d$ matrix, denoted~$\varrho$.
The map sends $x$ to the matrix
$\varrho \varrho^t$, which has rank at most~$d$.  (The rank of
$\varrho \varrho^t$ is less than~$d$
if $\prho$ has an affine span of dimension less than~$d$).
It is easy to check that $\varrho \varrho^t$ does not change if we
change $\varrho$ by an orthonormal transform, and that in fact the
entries of $\varrho \varrho^t$ are linear combinations of the squared
lengths of edges, so the map is well-defined and linear.
\end{proof}

\begin{definition}\label{def:faces}
  Let $K$ be a non-empty  convex set in an affine space.
  The
  \emph{dimension} of~$K$ is the dimension of the smallest affine
  subspace containing~$K$.  A point $x\in K$ is in the \emph{relative
    interior} $\Int(K)$ of~$K$ if there is a neighborhood $U$ of $x$ in 
  $\EE^{e}$
  so that $U\cap K\cong \RR^{\dim K}$.  (It is easy to see that the
  relative interior of $K$ is always non-empty.)

  The \emph{face}~$F(x)$ of a point $x\in K$ is the (unique) largest
  convex set contained in~$K$ so that $x \in \Int(F(x))$.  The faces
  of~$K$ are the sets $F(x)$ as $x$ ranges over~$K$.  (In
  Section~\ref{sec:extremity} we will
  see more characterizations of $F(x)$.)  Let $f(x)$ be the dimension
  of $F(x)$.
\end{definition}

\begin{lemma}
\label{lem:mes-faces}
Let $x = \ell_\Delta(\prho)$ with $\prho \in C^d(\Verts)$.
Then
$F(x) = \ell_\Delta(\Aff_d(\prho)) 
= \ell_\Delta(\GL_d(\prho))$.
When $\prho$ is proper,
$f(x) = \binom{d+1}{2}$.
\end{lemma}
Here $\GL_d$ is the group of $d$-dimensional linear transforms, and
$\Aff_d$ is the corresponding affine group.

\begin{proof}
By Lemma~\ref{lem:meas-sdp}, this reduces to understanding the faces
of the cone of $v-1$ by $v-1$
PSD matrices, which is well understood
(see, e.g.,~\cite[Theorem A.2]{pataki2000geometry}).
In particular, for $\varrho$ a $(v-1)\times d$ matrix,
the points in
$F(\varrho \varrho^t)$ in the semidefinite cone
are those matrices
of the form
$\varrho L^t L \varrho^t$ 
for some $d\times d$ matrix~$L$.  But each such
$\varrho L^t L \varrho^t$ 
maps under our
isometry  to 
$\ell_\Delta(\sigma)$ 
where   $\sigma \in C^d(\Verts)$ is obtained from $\prho$ by a $d$-dimensional
affine transform, with linear part described by $L$.
When $\prho$ has an affine span of dimension $d$,
$\varrho$ is of rank~$d$.  This implies that $F(\varrho \varrho^t)$ has
dimension $\binom{d+1}{2}$, as we can see by computing the dimension of
$\GL_d/O_d$, the linear transforms modulo the orthonormal transforms.
\end{proof}

\subsection{Dual vectors and PSD stresses}

\begin{definition}
Let us index each of the coordinates of 
  $\RR^{\binom{v}{2}}$ with an  integer pair
$ij$ with $1 \le i<j \le v$.
Given 
  $\phi$, a functional in the dual space
  $\bigl(\RR^{\binom{v}{2}}\bigr)^*$, define the  $v$-by-$v$
matrix $M(\phi)$ as follows: for $i \ne j$,
$M(\phi)_{ij}=
M(\phi)_{ji}\coloneqq  -\phi_{ij}$ and 
$M(\phi)_{ii} \coloneqq  -\sum_{j\neq i} M(\phi)_{ij}$
\end{definition}

\begin{lemma}
\label{lem:covec-stress}
  Let $\phi$ be a dual vector in $\bigl(\RR^{\binom{v}{2}}\bigr)^*$
 and let $\Omega\coloneqq M(\phi)$ be the
  corresponding 
matrix.
  Let $\prho \in C^d(\Verts)$.  Then
  \begin{equation}
    \langle \phi, \ell_\Delta(\prho) \rangle
= \frac{1}{2}\sum_{k=1}^d (\prho^k)^t \Omega \prho^k.
  \end{equation}
\end{lemma}

Here we use the notation $\prho^k$ for vector in $\RR^v$ 
describing the component of $\prho$ in the
$k$'th coordinate direction of $\EE^d$.

\begin{proof}
Both sides are equal to 
$\sum_k \sum_{i,j \mid i < j} (\prho^k_i-\prho^k_j)^2 \phi_{ij}$.
\end{proof}

\begin{definition} For $K\subset \RR^n$ a  convex set and
  $\phi$ a functional in the dual space
  $(\RR^n)^*$, we say that $\phi$ is \emph{tangent to 
$K$ at~$x \in K$} if
for all $y\in K$, $\langle \phi, y\rangle \geq 0$
while $\langle \phi, x\rangle =0$.
\end{definition}

\begin{lemma}
\label{lem:covec-psdstress}
  Let $\phi$ be a dual vector in $\bigl(\RR^{\binom{v}{2}}\bigr)^*$
that is tangent
to 
$M(\Delta_v)$ at $\ell_\Delta(\prho)$ for some $\prho \in C(\Verts)$.
Let $\Omega\coloneqq  M(\phi)$ be the
corresponding matrix.  
Then $\Omega$ is a PSD equilibrium stress matrix
for $(\prho,\Delta_v)$.
\end{lemma}
\begin{proof}
Conditions (1) and (3) in Definition~\ref{def:stress-matrix} are
automatic by definition of~$\Omega$, and condition~(2) is vacuous for
the complete graph $\Delta_v$.  It remains to check condition~(4) and
that $\Omega$ is PSD.

By Lemma~\ref{lem:covec-stress}, 
$\langle \phi, \cdot \rangle$  can be evaluated as a quadratic form
using the matrix $\Omega$.
Since $\langle \phi, \cdot \rangle$ 
is not negative anywhere on 
$M(\Delta_v)$, 
$\Omega$
must be PSD.
Because $\Omega$ is positive semidefinite
and
$\sum_{k=1} (\prho^k)^t \Omega \prho^k=0$,
it must also be true that for each $k$,
$(\prho^k)^t \Omega \prho^k = 0$ and so $\Omega \prho^k =0 $,
which is the last necessary condition to show $\Omega$ is an
equilibrium stress matrix.
\end{proof}

\begin{lemma}
\label{lem:covec-rankpsdstress}
In the setting of Lemma~\ref{lem:covec-psdstress}, suppose
furthermore that
$\phi$ is  only  tangent to  points in $F(\ell_\Delta(\prho))$,
and suppose the affine span of $\prho$ has dimension $d$.

Then $\Omega$ is a PSD equilibrium stress matrix
for $(\prho,\Delta_v)$ with rank $v-d-1$.
\end{lemma}
\begin{proof}
From Lemma~\ref{lem:covec-psdstress}, $\Omega$ is a 
PSD equilibrium 
stress matrix for $(\prho,\Delta_v)$. 
By assumption, its kernel is 
spanned by  the coordinates
of frameworks that map under $\ell_\Delta$ to $F(\ell_\Delta(\prho))$.
From Lemma~\ref{lem:mes-faces}, such frameworks 
are $d$-dimensional affine
transforms of $\prho$. Thus the kernel is spanned by the 
$d$ coordinates of $\prho$ and the all-ones vector.  Since the kernel
has dimension $d+1$, $\Omega$ has rank $v-d-1$.
\end{proof}

\begin{lemma}
\label{lem:covec-gamma-rankpsdstress}
In the setting of Lemma~\ref{lem:covec-rankpsdstress}, suppose
furthermore that
$\phi_{ij}=0$ for $\{i,j\} \not\in \Edges(\Gamma)$ for a graph~$\Gamma$.

Then $\Omega$ is a positive semidefinite equilibrium stress matrix
for $(\prho,\Gamma)$ with rank $v-d-1$.
\end{lemma}
\begin{proof}

$\Omega$ must have zeros in all coordinates corresponding to non-edges
in $\Gamma$, thus it will be 
 an equilibrium  stress matrix for $\Gamma$ as well
as for $\Delta_v$.
The rest follows from Lemma~\ref{lem:covec-rankpsdstress}.
\end{proof}

To be able to use  Lemma~\ref{lem:covec-gamma-rankpsdstress}, we want
to find a 
dual vector $\phi$ that is  only  tangent to  points in 
$F(\ell_\Delta(\prho))\subset M(\Delta_v)$.
Additionally we need 
that $\phi_{ij}=0$ for 
$\{i,j\} \not\in \Edges$. 
As we will see, this will correspond to finding 
a dual vector in $(\RR^e)^*$ that is only tangent to 
points in $F(\ell_\Gamma(\prho)) \subset M(\Gamma)$.
The proper language for this kind of condition is the notion of convex
extreme and exposed points, which we turn to next.


\section{Convex sets and projection}
\label{sec:convex}

We will prove our main theorems in the more general setting of closed convex
cones.  Our goal is the following general theorem about projections of
convex sets, which will quickly imply
Theorem~\ref{thm:generic-ur-min-stress-kernel}.  (See
Definitions~\ref{def:extreme}, \ref{def:exposed},
and~\ref{def:univ-rigid} for the terms used.)

\begin{theorem}
  \label{thm:urexp}
  Let $K$ be a closed
  line-free convex semi-algebraic set in $\RR^m$, and $\pi: \RR^m
  \to \RR^n$ a projection, both defined over~$\QQ$.
  Suppose $x$ is locally generic in $\ext_k(K)$ 
  and
  universally rigid under $\pi$.
  Then $\pi(x)$ is $k$-exposed.
\end{theorem}

\subsection{Extreme points}
\label{sec:extremity}

\begin{definition}\label{def:extreme}
Let $K$ be a non-empty, convex set.
  A point $x\in K$ is
  \emph{$k$-extreme} if $f(x) \le k$.  (Recall from
  Definition~\ref{def:faces} that $f(x)$ is the dimension of the face
  of $K$ containing~$x$.)
Let $\ext_k(K)$ be the set of
  $k$-extreme points in~$K$.  
\end{definition}

We will also use the following elementary propositions
(see for example 
the exercises in chapter 2.4 of~\cite{grunbaum2003convex}.)

\begin{proposition}
  For $K$ a convex set and $x \in K$, the following statements
  are equivalent:
  \begin{itemize}
  \item $f(x) \le k$, and
  \item $x$ is not in the relative interior of any non-degenerate
    $(k+1)$-simplex contained in~$K$.
  \end{itemize}
\end{proposition}

\begin{proposition}\label{prop:face-segment}
\label{prop:altFace}
  For $K$ a convex set and $x\in K$, the face $F(x)$ is the set
  of points $z \in K$ so that there is an $z' \in K$ with $x$ in the
  relative interior of the segment $[z', z]$.
\end{proposition}

\begin{remark}
  One special case of Proposition~\ref{prop:face-segment} is when
  $z=x$, in which case we also take $z'=x$ and the segment consists of
  a single point.
\end{remark}

\begin{corollary}
  For $K$ a convex set, the faces of~$K$ are the  convex
  subsets~$F$ of~$K$ such that every line segment in $K$ with a relative
  interior point in $F$ has both endpoints in~$F$.
\end{corollary}

\subsection{Exposed points}
\label{sec:exposedness}

\begin{definition}\label{def:exposed}
  A point $x\in K$ is
\emph{$k$-exposed} if there
  is a closed half-space~$H$ containing~$x$ so that $\dim H \cap K \le
  k$.  Let $\exp_k(K)$ be the set of $k$-exposed points in~$K$.
\end{definition}

\begin{figure}
  \begin{center}
    $\mfigb{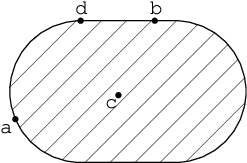}$
  \end{center}
  \caption{Examples of extreme and exposed points.  For the
    racetrack shape illustrated (bounded by semi-circles and line
    segments), $a$ is $0$-extreme and $0$-exposed, $b$ is
    $1$-extreme and $1$-exposed, $c$ is $2$-extreme and $2$-exposed,
    and $d$ is $0$-extreme and $1$-exposed.}
  \label{fig:examples-exposed}
\end{figure}

See Figure~\ref{fig:examples-exposed} for some examples of $k$-extreme
and $k$-exposed points, including a case where they differ.
If $x$ is $k$-exposed, it is also $k$-extreme, as any simplex
containing $x$
in its relative interior is contained in any supporting hyperplane.

The following theorem is a crucial tool in our proof.

\begin{otheorem}[Asplund~\cite{Asplund63:k-extreme}]
  \label{thm:k-exposed-dense}
  For $K$ a closed convex set, 
the $k$-exposed points are dense in the $k$-extreme points.
\end{otheorem}

\begin{remark}
  Asplund only states and proves Theorem~\ref{thm:k-exposed-dense} for
  compact convex sets.  The result follows for any closed 
 convex
  set~$K$ as follows.  For $x\in \ext_k(K)$, fix $R > 0$ and consider
  the compact convex set $K' \coloneqq K \cap \overline{B}_R(x)$.  
Then Asplund's
  theorem says we can find $k$-exposed
  points~$z$ in~$K'$ that are arbitrarily close to~$x$.  If $z$ is
  sufficiently close to~$x$ (say, within $R/2$ of $x$), it will also be
  $k$-exposed in~$K$, as desired.
\end{remark}

Theorem~\ref{thm:k-exposed-dense} is an extension of a theorem of
Straszewicz~\cite{Straszewicz35:exposed-dense}, who proved it in the
case $k=0$.  There are several improvements of
Theorem~\ref{thm:k-exposed-dense}, replacing ``dense'' with
various stronger assertions. In our
applications we will consider semi-algebraic 
sets~$K$ defined
over~$\QQ$.  Note that for any such set~$K$, $\ext_k(K)$ and
$\exp_k(K)$ are also
semi-algebraic over~$\QQ$, as the $k$-extreme and $k$-exposed
conditions can be phrased
algebraically.

\begin{corollary}\label{cor:generic-extreme-exposed}
Let  $K$ be  a  closed, convex,
semi-algebraic
set.
If $x$ is
  locally generic in $\ext_k(K)$, then $x$ is
  $k$-exposed.
\end{corollary}

\begin{proof}
This follows from an application of Lemma~\ref{lem:dense}
 to the inclusion $\exp_k(K)
  \subset \ext_k(K)$ in a small neighborhood around $x$.
\end{proof}

\subsection{Projection}
\label{sec:projection}

The two convex sets $M(\Delta_v)$ and
$M(\Gamma)$ are related by a projection from $\RR^{\binom{n}{2}}$
to~$\RR^e$.  We will consider this
situation more generally.

Throughout this section, fix two Euclidean spaces,
$\EE^m$ and $\EE^n$, with a surjective projection map $\pi: \EE^m \to
\EE^n$ and a closed convex set $K \subset \EE^m$.  
We will work with the image $\pi(K)
\subset \EE^n$.
Let $\pi_K$ be the
restriction of $\pi$ to $K$.

\begin{remark}
  In general, $\pi(K)$ need not be closed, even if $K$ is.  We
  continue to work with $\pi(K)$ and its ``faces'' as defined in
  Definition~\ref{def:faces}, even if it is not
  closed.
\end{remark}

\begin{definition}\label{def:univ-rigid}
  We say that $x\in K$ is \emph{universally rigid (UR) under
    $\pi$} if $\pi_K^{-1}(\pi(x))$ is
  the single point~$x$.
\end{definition}

Note that for $y \in \pi(K)$, $\pi_K^{-1}(y)$ is always a convex set
in its own right.

\begin{lemma}\label{lem:faces-nest}
Let $K$ be a convex set and $\pi$ a projection.
  For any $x\in K$, $\pi(F(x)) \subset F(\pi(x))$.
\end{lemma}
\begin{proof}
  The set $\pi(F(x))$ is a convex set containing~$\pi(x)$ in its relative
  interior, so by definition is contained in $F(\pi(x))$.
\end{proof}

A point that is $0$-extreme in a convex set is called a 
\emph{vertex} of the set.
A convex set is \emph{line free} if it contains no
complete affine line.
Recall that a non-empty, closed, line-free convex set has a vertex
(see e.g.,~\cite[2.4.6]{grunbaum2003convex}).
Our convex
sets of interest, $M(K_\Delta)$ and $M(\Gamma)$, are closed and also
automatically
line free, as the length-squared takes only positive values and any
line contains both positive and negative values in at least one coordinate.

\begin{proposition}
  \label{prop:dims-decrease}
Let $K$ be a convex set and $\pi$ a projection.
  For any $y \in \pi(K)$ and  $x$  a vertex of 
$\pi_K^{-1}(y)$, $F(x)$ maps
  injectively under $\pi_K$ into $F(y)$.  In particular, $f(x) \le
  f(y)$. If $K$ is closed and line free, then for every $y \in \pi(K)$
  there is a $x \in \pi^{-1}(y)$ so that $f(x) \le f(y)$.
\end{proposition}

\begin{proof}
Let $A(x)$ be
the smallest affine subspace containing $F(x)$. 
Suppose $A(x)$
contained a  direction vector $v$
in the kernel of $\pi$. Then, since $x$ is in the relative interior
of $F(x)$ for small enough $\epsilon$, we would have the segment
$[x+\epsilon v,x-\epsilon v]$ contained in $F(x)$ and also in 
$\pi_K^{-1}(y)$.
This would contradict
the assumption that 
 $x$ is  a vertex of 
$\pi_K^{-1}(y)$.
Thus 
  $F(x)$ maps injectively under~$\pi_K$.  

From
  Lemma~\ref{lem:faces-nest} we see that $f(x) \le f(y)$.  For the
  last part of the lemma, observe
  that if $K$ is closed and
  line free, so is $\pi_K^{-1}(y)$, and so $\pi_K^{-1}(y)$
  has a vertex, which is the desired point~$x$ by the first part.
\end{proof}

\begin{proposition}
\label{prop:dims-increase}
Let $K$ be a convex set and $\pi$ a projection.
  For any $y \in \pi(K)$ and  $x$  
in the relative
interior of $\pi_K^{-1}(y)$,
$F(x) =  \pi_K^{-1}(F(y))$.  
In particular,
$f(x) \ge f(y)$.
\end{proposition}
\begin{proof}
  By Lemma~\ref{lem:faces-nest}, we already know that $F(x) \subset
  \pi_K^{-1}(F(y))$, so we just need to show the other inclusion.
  Pick any $y_1 \in F(y)$ and $x_1 \in
  \pi_K^{-1}(y_1)$.  We must show that $x_1 \in F(x)$.  Since 
$y_1
  \in F(y)$, 
by Proposition~\ref{prop:face-segment} there is an
  $y_2 \in \pi(K)$ so that $y \in \Int([ y_1,y_2])$.
  Let $x_2 \in \pi_K^{-1}(y_2)$.
  Let $x' \in K$ be the
  (unique) point of intersection of $[x_1, x_2]$ with $\pi_K^{-1}(y)$.
  Since $x \in \Int(\pi_K^{-1}(y))$, there is a point
  $x''\in\pi_K^{-1}(y)$ with $x \in \Int([x', x''])$.  But then
  $x$ is in the interior of the simplex $[x_1, x_2, x'']$, which implies
  that $x_1\in F(x)$.
\end{proof}

In order to be able to apply Asplund's theorem at $\pi(x)$ we
need $\pi(x)$ to be locally generic in $\ext_k(\pi(K))$. The following
lemma will assist us.

\begin{lemma}
  \label{lem:proper-map-inverse}
  Let $K$ be a closed, convex set in $\RR^m$
and let
  $\pi_K: K \to \RR^n$ be a projection.  
Let $x\in K$ be universally rigid under $\pi$. 
Then for any $\epsilon$  there is a $\delta > 0$ so that 
$\pi_K^{-1}(B_\delta(\pi(x))) \subset
 B_\epsilon(x) \cap K$.
\end{lemma}

Intuitively, the only points in~$K$ that map close to $\pi(x)$ are
close to~$x$.

\begin{proof}
Suppose not. Then
there is an $\epsilon$ and 
a sequence of points $x_i$ with the following property:
$\pi(x_i)$ approach $\pi(x)$, while 
$d(x,x_i) > \epsilon$. 
Let $x_i'$ be the point on the interval $[x,x_i]$ that is
a distance exactly $\epsilon$ from $x$. By convexity $x_i' \in K$.
Then the $x_i'$ are in a bounded and closed set and therefore have
an accumulation point $x'\in K$, which will also have distance
$\epsilon$ from~$x$. 
By the linearity of $\pi$
the $\pi(x'_i)$ also approach $\pi(x)$.
So by continuity $\pi(x')=\pi(x)$. 
This contradicts the universal rigidity of~$x$.
\end{proof}

And now we can prove the following.

\begin{lemma}
  \label{lem:loc-gen}
  Let $K$ be a closed line-free
  convex semi-algebraic set in $\RR^m$ and $\pi: \RR^m
  \to \RR^n$ a projection, both defined over~$\QQ$.
  Suppose $x$ is locally generic in $\ext_k(K)$ 
  and
  universally rigid under $\pi$.
  Then $\pi(x)$ is locally generic in $\ext_k(\pi(K))$.
\end{lemma}
\begin{proof}
  By Proposition~\ref{prop:dims-decrease}, 
  $\ext_k(\pi(K)) \subset \pi(\ext_k(K))$.
  By Proposition~\ref{prop:dims-increase} 
and the universal rigidity
  of $x$, we have that $\pi(x)$ is in $\ext_k(\pi(K))$.

  Let $V\coloneqq  \ext_k(K) \cap B_\epsilon(x)$ and 
  $W\coloneqq  \ext_k(\pi(K)) \cap B_\delta(\pi(x))$.
For 
  sufficiently small~$\epsilon$,
 by local genericity,  $x$ is generic in $V$.
Meanwhile, $W \subset \pi(\ext_k(K))$, and,
from Lemma~\ref{lem:proper-map-inverse} ,
for  small enough $\delta$, we have
$\pi_K^{-1}(W) \subset B_\epsilon(x)$ and thus
$W \subset \pi(V)$. 

Thus from Lemma~\ref{lem:image-generic}, 
we have  $\pi(x)$ generic
in $W$. Thus $\pi(x)$ is locally generic in 
$\ext_k(\pi(K))$.
\end{proof}

\begin{remark}
In our special case where $K$ is the cone of positive 
semidefinite matrices, $x$ is in fact generic in 
$\ext_k(K)$ (which is irreducible), and thus
the second paragraph in the above proof (and in turn
Lemma~\ref{lem:proper-map-inverse}) are not needed.
\end{remark}

Finally, in order to be able to apply Asplund's theorem to $\pi(K)$ we
need $\pi(K)$ to be a closed set. In the graph embedding case,
we can use the fact that $\pi$ is a proper map whenever $\Gamma$ is 
connected,
and thus $\pi(K)$ must
be closed. 
In the setting 
of a general $K$ and $\pi$
we can argue closedness using standard techniques.

\begin{definition}
A direction of recession $v$ for a convex set~$K$ is a vector such that 
for every (equivalently, any) $x \in K$ and every $\lambda \geq 0$ we have
$x + \lambda v \in K$ (see e.g.~\cite[Chapter 8]{rockefeller1970convex}).
\end{definition}

A sufficient condition for closedness
is given by the following 
theorem~\cite[Theorem 9.1]{rockefeller1970convex}.
\begin{otheorem}
\label{thm:fibcl}
Let $K$ be a line free closed convex set and $\pi$ a projection.
If $K$ does not have a direction of recession $v$ with
$\pi(v)=0$, then 
then $\pi(K)$ is closed.
\end{otheorem}

\begin{corollary}
\label{cor:projclosed}
Let $K$ be a line free closed convex set and $\pi$ a projection.
If there is a point $x\in K$ that is universally rigid under $\pi$,
then $\pi(K)$ is closed.
\end{corollary}
\begin{proof}

Suppose there was a direction of recession $v$ for $K$ with
$\pi(v)=0$. 
Then $\pi(x + v) = \pi(v)$ with $x+v$ in $K$. This would contradict
the universal rigidity of $x$.
\end{proof}

Putting this all together, we can deduce Theorem~\ref{thm:urexp}.

\begin{proof}[Proof of Theorem~\ref{thm:urexp}]
  From Corollary~\ref{cor:projclosed}, $\pi(K)$ is closed.
  From Lemma~\ref{lem:loc-gen}, $\pi(x)$ is 
locally generic in $\ext_k(\pi(K))$.
Thus by
  Corollary~\ref{cor:generic-extreme-exposed}, $\pi(x)$ is $k$-exposed.
\end{proof}

\subsection{Proof of main theorem}
We are now in position to complete the proof of
Theorem~\ref{thm:generic-ur-min-stress-kernel}.
Recall that we have a universally rigid framework $(p,\Gamma)$ with
$p$ generic in $C^d(\Verts)$.
Let $K$ be  $M(\Delta_v)$, 
and $\pi$ be its projection to 
$M(\Gamma)$, $k$ be $\binom{d+1}{2}$.
and $x \in M_d(\Delta_v) \subset K$ be $\ell_\Delta(\prho)$.
The assumption that $\prho$ is generic in $C^d(\Verts)$ implies that
$x$ is generic in $M_d(\Delta_v) = \ext_k(K)$.

From Theorem~\ref{thm:urexp}, 
$\pi(x)$ is $k$-exposed. 
So
there must be a closed halfspace $H$ in $\RR^e$ whose 
intersection with $\pi(K)$
is exactly $F(\pi(x))$. 

The preimage $\pi^{-1}(H)$
is
a closed halfspace  
in $\RR^{\binom{v}{2}}$
whose intersection with $K$
is exactly $\pi_K^{-1}(F(\pi(x)))$, which
by Proposition~\ref{prop:dims-increase} is
$F(x)$. 
Since $K$ is a cone, 
$\pi^{-1}(H)$ must also include  the origin.
When  $v>d+1$, then $F(x)$ is not the entire 
measurement cone $M(\Delta_v)$,
and the halfspace 
$\pi^{-1}(H)$ must not include any interior points of $M(\Delta_v)$.

Thus 
$\pi^{-1}(H)$ is represented by a dual vector 
$\phi\in \bigl(\RR^{\binom{v}{2}}\bigr)^*$, in the sense that
$\pi^{-1}(H) = \bigl\{\,x \in \RR^{\binom{v}{2}} \mathbin{\big\vert} \phi(x) 
\leq  0\,\bigr\}$.
Moreover, this $\phi$ is tangent to  $M(\Delta_v)$ exactly at $F(x)$.

Since we started with a halfspace in $\RR^e$,
$\phi_{ij} =0$ for
$\{i,j\} \not\in \Edges$.
From Lemma ~\ref{lem:covec-gamma-rankpsdstress}, then,
$\Omega\coloneqq  M(\phi)$ must be  a positive semi-definite 
equilibrium stress matrix of $(\prho,\Gamma)$
with rank $v-d-1$.


\section{Relation to semidefinite programming}
\label{sec:sdp}

Theorem~\ref{thm:generic-ur-min-stress-kernel}
can be interpreted as a strict complementarity
statement about a particular
semidefinite program (SDP) and its associated dual program, as we
will now explain.
For a general survey on semidefinite programming 
see Vandenberghe and Boyd~\cite{vandenberghe1996semidefinite}. Our notation is 
similar to that of Pataki and Tun{\c{c}}el~\cite{pataki2001generic}.
For a related discussion on universal rigidity, semidefinite programming
and complementarity, see, for instance, the paper by So and
Ye~\cite{SY07:SemidefProgramming}.

\subsection{Semidefinite programming}
\label{sec:sdp-graphs}

We first recall the basic definitions of semidefinite programming.
Let $S^n$ be the linear space of all symmetric $n$-by-$n$ matrices and
$S^n_+ \subset S^n$ be the cone of
positive semidefinite matrices.
A \emph{semidefinite program} is given by a triple $(L,b,\beta)$ where
$L \subset S^n$ is a linear subspace, $b \in S^n$ is a matrix, and
$\beta \in (S^n)^*$ is a dual to a matrix (which we can also view as a
matrix, using the natural inner product $\langle A, B \rangle = \tr
(A^t B) = \sum_{i,j} A_{ij} B_{ij}$).
This describes a primal constrained optimization problem:
\[
\inf_x \{\langle \beta,x\rangle \mid x \in (L+b) \cap S^n_+
\}
\]
That is, we optimize a linear functional over all symmetric positive
semidefinite matrices~$x$, subject to the constraint
that $x$ lies in a chosen affine space.

Associated with every primal SDP  is the
(Lagrangian) dual program
\[
\inf_\Omega \{\langle  \Omega,b\rangle \mid \Omega \in (L^{\perp}+\beta) 
\cap (S^n_+)^*\}
\]
Since the cone of PSD matrices is self-dual, elements of $(S^n_+)^*$
also correspond to PSD matrices, so this dual program can be thought
of as a semidefinite program itself.

\begin{definition}
We say a pair of 
primal/dual points 
$(x,\Omega)$ with $x \in S^n_+$ and $\Omega \in (S^n_+)^*$
are a \emph{complementary pair} for the program $(L,b,\beta)$
if both are feasible points for the respective programs and
$\langle \Omega,x\rangle=0$.
\end{definition}

A calculation shows that 
for any feasible primal/dual pair 
of points
$(x,\Omega)$, we have 
$\langle \beta,b\rangle -\langle  \Omega,b\rangle \leq
\langle \beta,x\rangle $
and that complementarity implies 
$\langle \beta,b\rangle -\langle  \Omega,b\rangle =
\langle \beta,x\rangle$.

Thus, writing 
the dual problem in the alternative form
\[
\sup_\Omega \{
\langle \beta,b\rangle -\langle  \Omega,b\rangle 
\mid \Omega \in (L^{\perp}+\beta) 
\cap (S^n_+)^*\}
\] 
we see that a 
 complementary pair must represent optimal solutions to 
the primal and dual SDPs respectively.

\begin{definition}
We say that an SDP problem $(L,b,\beta)$ is \emph{gap free} if
it has complementary pairs of solutions.
\end{definition}

Gap freedom is typically seen as a fairly mild constraint on an SDP
problem.

For every  complementary pair  $(x,\Omega)$ 
we have~\cite{alizadeh1997complementarity}
\begin{equation}
\rank(x) + \rank(\Omega) \leq n.\label{eq:complement-rank-bound}
\end{equation}

\begin{definition}
A complementary pair  $(x,\Omega)$ is said to be \emph{strictly
  complementary} if we have
$\rank(x) + \rank(\Omega)= n$. 
\end{definition}

There is a rich theory on when SDPs have strictly complementary
solution pairs.  (See \cite{alizadeh1997complementarity,pataki2001generic,
  nie2010algebraic}; see also Section~\ref{thm:AHOsc}.)
In particular, certain SDP algorithms converge
more quickly on problems that satisfy strict complementarity~\cite{luo1998superlinear}.  
The linear programming counterpart of strict complementarity is strict
complementary slackness~\cite{dantzig2003linear}.
In the linear case, strict complementarity
is always achieved.  By contrast, for SDP problems it depends on the
particular problem.

\subsection{Graph embedding as an SDP}

In the context of graph embedding, suppose
we are looking for a 
configuration of unconstrained dimension, $\prho \in C^{v}(\Verts)/\Eucl(v)$,
such that the framework of $\Gamma$
is constrained to have 
a set of squared edge lengths 
$d^2_{ij}$ for $ij\in\Edges(\Gamma)$.
It is well known that we can 
set up this graph embedding problem 
as a semidefinite program~\cite{linial1995gga}, as we will now review.

To a configuration~$p$ we can associate a $v \times v$ \emph{Gram
matrix}~$x$, with $x_{ij} = \prho(i) \cdot \prho(j)$.  The Gram
matrix of~$p$ is unchanged by elements of $O(v)$ (although it does change
when $p$ is translated).
The matrix rank, $r$, of such an 
$x$ is the dimension of the 
linear span of the associated configuration~$\prho$.
This will be typically be $d+1$; 
one greater than the affine span, $d$,  of the framework.

In the graph embedding SDP,
we set $n\coloneqq v$ and take the Gram matrix $x$ as our unknown.
The distance constraint at an edge
$ij$ can be expressed as the linear constraint
$x_{ii}+x_{jj}-2x_{ij}=(d_{ij})^2$. The collection of matrices satisfying
these constraints for all edges forms an affine space; we choose $L$
and $b$ such  that
$L+b$ is this affine space. 
In our context, we are only interested
in feasibility, and thus have
no objective function, so we set $\beta \coloneqq  0$. 
Note that semidefinite programs do
not allow us to explicitly
constrain the rank of
the solution (which is related to the dimension of the configuration).

Let us now look at the dual program to our
graph embedding SDP problem.
In the primal problem, the linear space $L$ corresponds to 
symmetric matrices $x$ 
with
$x_{ii}+x_{jj}-2x_{ij}=0$ for all $\{i,j\} \in \Edges(\Gamma)$. 
Thus the space $L^\perp$,
when represented as matrices, is spanned by a basis elements $B_{ij} =
-e_{ii}+e_{ij}+e_{ji}-e_{jj}$ for $\{i,j\}\in\Edges(\Gamma)$, where
$e_{ij}$ is the elementary matrix with a $1$ in the $ij$ entry and $0$ elsewhere.
The space $L^\perp$ is therefore those
matrices $\Omega$
with row-sums of zero and with
$\Omega_{ij} = 0$ for all $i \ne j$, $\{i,j\} \notin \Edges(\Gamma)$.
That is, $L^\perp$ is the
space of all (not necessarily equilibrium) stress matrices.
We also have $\beta=0$, so we optimize over these matrices.

Minimizing $\langle \Omega, b\rangle$ imposes
the   further constraint that the solution be an equilibrium
stress matrix.
This can be seen as follows.
Since $\beta=0$, the all-zeros matrix  is dual feasible.
Whenever there is some primal feasible $x$ corresponding to a
configuration~$p$, the SDP must be gap free
as 
$\langle 0, x\rangle=0$, and the optimal $\Omega$ are the
feasible $\Omega$ with
$\langle \Omega, x\rangle=0$.
For any $\Psi \in (S^n_+)^*$
we have
$\langle \Psi, x\rangle= 
\frac{1}{2}\sum_{k=1}^r (\prho^k)^t \Psi \prho^k$.
(Here we use the notation $\prho^k$ for vector in $\RR^n$ 
describing the component of $\prho$ in the
$k$'th coordinate direction of $\RR^r$.)
And thus, as in Lemma~\ref{lem:covec-psdstress},
the feasible $\Omega$ with 
 $\langle \Omega, x\rangle=0$ are 
the PSD 
equilibrium stress matrices for $(\prho,\Gamma)$.

In this language we can restate
Theorem~\ref{thm:generic-ur-min-stress-kernel} as follows.

\begin{proposition}
Let $\prho$ be generic in $C^d(\Verts)$,
let $(\prho,\Gamma)$ be 
universally rigid in $E^d$, and let $(L,b,0)$ be the associated SDP
using the graph $\Gamma$ and the distances of the
edges in $(\prho,\Gamma)$.
Then $(L,b,0)$ has
a strictly complementary solution pair $(x,\Omega)$.
\end{proposition}
\begin{proof}
Simply let $x$  be the 
Gram matrix 
corresponding to a
 translation of $\prho$ such that its affine span does not include
the origin.  Thus $x$ will have
rank $d+1$. Let $\Omega$   be 
a PSD equilibrium stress matrix of rank $v-d-1$, which exists
by Theorem~\ref{thm:generic-ur-min-stress-kernel}.  Then $\rank x +
\rank \Omega = v$.
\end{proof}

\subsection{SDP feasibility}

Our discussion of universal rigidity and 
equilibrium stress matrices carries
over directly to any SDP feasibility problem $(L,b,0)$ (i.e., where $\beta=0$).
To see this we first set up some notation.  Pick an integer $r > 0$ and
let $k\coloneqq  \binom{r+1}{2}$.
Then, as in Lemma~\ref{lem:mes-faces}, $\ext_k(S^n_+)$ 
is the set of the PSD matrices of rank less than or
equal to $r$, which we will denote $S^n_{+r}$.
Let $\pi$ be projection from $S^n$ to $S^n/L$, and $\pi_+$ its restriction
to $S^n_+$.
A point $x \in S^n_+$ is a solution to the feasibility problem
iff 
$\pi(x) = \pi(b)$, or $x \in \pi_+^{-1}(\pi(b))$.

We say that a point $x \in S^n_+$ is universally rigid under $\pi$
iff $\pi_+^{-1}(\pi(x))$ is a single point: A universally rigid
feasible point is the unique feasible solution to
$(L,b,0)$.  As we vary~$b$ with fixed~$L$, we will see as solutions all points
in~$S^n_+$; it may happen that for some~$b$ there is a unique, generic
feasible solution.

\begin{proposition}
\label{pro:sdp}
Let $(L,b,0)$ be an SDP feasibility problem, where $L$ has rational
coefficients.  Suppose there is a unique feasible solution~$x$
which is
generic in $S^n_{+r}$.
Then there exists an optimal dual solution $\Omega$ such that
$(x,\Omega)$ is a 
strictly complementary pair.
\end{proposition}
\begin{proof}
Since $x$ is generic in $S^n_{+r} = \ext_k(S^n_+)$,
by Theorem~\ref{thm:urexp} $\pi(x)$ is 
$k$-exposed.
Thus there is a PSD matrix
$\Omega$ in 
$(S^n/L)^*=L^\perp$ 
that is tangent to $\pi(S^n_+)$
at $\pi(x)$, with contact $F(\pi(x))$. 
Since $\Omega  \in L^\perp$, it is a feasible point of the dual SDP\@.
Since 
$\Omega$ 
is tangent to $\pi(S^n_+)$
at $\pi(x)$, it is also tangent 
to $S^n_+$
at $x$;
thus $\langle \Omega, x\rangle=0$ and $\Omega$ is
a complementary and therefore optimal dual solution.

The tangency of $\Omega$ 
to  $S^n_+$
at $x$ has contact
$\pi^{-1}_+(F(\pi(x))$ which,
by 
Proposition~\ref{prop:dims-increase}, 
is $F(x)$.
Let $\prho$ be an $r$-dimensional configuration of $n$ points
with Gram matrix $x$.  From the facial stucture of the PSD cone
we see that $F(x)$ consists only of 
Gram matrices of $r$-dimensional linear transforms of~$\prho$.
We also have the relation 
$\langle \Omega, x\rangle= 
\frac{1}{2}\sum_{k=1}^r (\prho^k)^t \Omega \prho^k = 0$.
Thus 
we see that $\Omega$ must have 
a kernel of dimension $r$ and 
have rank 
$n-r$.
\end{proof}

\subsection{SDP optimization}

We can apply the approach above to a more general SDP problem
$(L,b,\beta)$, where $\beta \neq 0$.  Specifically, we
prove the following.

\begin{theorem}
\label{thm:sdp-opt}
Let $(L,b,\beta)$ be a gap-free SDP problem, with $L$ and $\beta$ rational.
Suppose there is a unique optimal solution $x$ which is generic
in $S^n_{+r}$.
Then there exists a strictly complementary pair $(x, \Omega)$.
\end{theorem}

We will prove the theorem by reducing to the $\beta=0$ case.  For any
SDP optimization problem $(L,b,\beta)$ with an optimal solution~$x$, there is an
associated SDP feasibility problem $(L', b', 0)$, where $L' + b'$ is
$\{ y \in L + b \mid \langle \beta, y\rangle = \langle \beta, x\rangle\}$.
In particular, $L'$ is $L \cap \ker \beta$ (which is still rational) and the dual space $L'^\perp$ is  $L^\perp +
\langle\beta\rangle$.
Then feasible solutions to $(L', b', 0)$ correspond to optimal
solutions to $(L, b, \beta)$.

\begin{lemma}
\label{lem:sdp-almostSC}
Let $(L,b,\beta)$ be an SDP problem, with $L$ and $\beta$ rational.
Suppose there is a unique optimal solution~$x$ which is generic 
in $S^n_{+r}$.
Then there exists a PSD $\Omega$ such that 
$\rank \Omega = n-r$, 
with 
$\Omega \in L^\perp + \langle\beta\rangle$
and 
$\langle \Omega,x\rangle=0$.
\end{lemma}
\begin{proof}
Apply Proposition~\ref{pro:sdp} to the feasibility problem $(L', b',
0)$ constructed above.  This gives a $\Omega \in L'^\perp = L^\perp +
\langle \beta \rangle$ that is strictly complementary to~$x$.
\end{proof}

Note that we have not yet proved Theorem~\ref{thm:sdp-opt}, as
Lemma~\ref{lem:sdp-almostSC} gives $\Omega$ in the linear space
$L^\perp + \langle\beta\rangle$ rather than the desired affine space
$L^\perp + \beta$.

\begin{proof}[Proof of Theorem~\ref{thm:sdp-opt}]
Let $\Omega_1$ be the PSD matrix given by
Lemma~\ref{lem:sdp-almostSC}.
From the assumption of gap freedom
there exists a PSD   $\Omega_2$ 
with 
$\Omega_2 \in L^\perp + \beta$
and $\langle \Omega_2,x\rangle=0$.
Thus 
for any positive scalars $\lambda_1$ and $\lambda_2$,
the matrix $\Omega \coloneqq  \lambda_1 \Omega_1 + \lambda_2 \Omega_2$
is PSD,
has $\rank$ no less than $n-r$, 
and has $\langle \Omega,x\rangle=0$.
By adjusting $\lambda_1$ and $\lambda_2$ we can achieve
$\Omega \in L^\perp +  \beta$.
By Equation~\eqref{eq:complement-rank-bound}, we in fact have $\rank
\Omega = n-r$.
\end{proof}

We note though that for a given $(L,\beta)$ and rank $r$, 
as we vary $b$, 
there
may be no 
unique solutions with rank $r$.
Even if there are such 
$x$, they may all  be non-generic. 
In such cases, Theorem~\ref{thm:sdp-opt} is vacuous, as there are no
choices of $(L, b, \beta)$
satisfying the hypotheses.

\subsection{Genericity}
\label{sec:genericity}

In the context of SDP optimization, the hypothesis in
Theorem~\ref{thm:sdp-opt} that $L$ and $\beta$ be rational may seem a
little unnatural.  This hypothesis can be relaxed if we
work with points that are generic over a field that is larger
than~$\QQ$.

\begin{definition}
  Let $\kk$ be a field containing $\QQ$ and contained in $\RR$.  A
  semi-algebraic set~$S \subset \RR^n$ is \emph{defined over $\kk$} if
  there is a set of equalities and inequalities defining~$S$ with
  coefficients in~$\kk$.  If $S$ is defined over~$\kk$, a point $x \in
  S$ is \emph{generic over~$\kk$} if the coordinates of~$x$ do not
  satisfy any algebraic equation with coefficients in $\kk$ beyond
  those that are satisfied by every point in~$S$.  Similarly, $x$ is
  \emph{locally generic over~$\kk$} if for small enough $\epsilon$,
  $x$ is generic in $S \cap B_\epsilon(x)$.

  A \emph{defining field} of a semi-algebraic set~$S$, written
  $\QQ[S]$, is any field over which it is defined.  (There is a unique
  smallest field for algebraic sets by a result of Weil
  \cite[Corollary IV.3]{Weil46:FoundAlgGeom}.  For our purposes, we can use any
  field over
  which $S$ is defined.)
  Similarly, if $f : X \to Y$ is a map between semi-algebraic sets,
  then $\QQ[f]$ is a field over which it is defined (or,
  equivalently, a field over which the graph of~$f$ is
  defined).
\end{definition}

For instance, if $S$ is a single point~$x$, we can take $\QQ[S]$ to be
the same as
$\QQ[x]$, the smallest field containing all the coordinates of~$x$.
Also, if $x$ is generic over $\QQ$, then $y$ is generic over $\QQ[x]$
iff the pair $(x,y)$ is generic over~$\QQ$.

With this definition, Theorems \ref{thm:urexp} and~\ref{thm:sdp-opt}
can be improved to allow non-rational sets and projections.

\begin{citingthm}[\ref*{thm:urexp}$'$]
  Let $K$ be a closed
  line-free convex semi-algebraic set in $\RR^m$, and $\pi: \RR^m
  \to \RR^n$ a projection.
  Suppose $x$ is locally generic over $\QQ[K,\pi]$ in $\ext_k(K)$ 
  and
  universally rigid under $\pi$.
  Then $\pi(x)$ is $k$-exposed.
\end{citingthm}

\begin{citingthm}[\ref*{thm:sdp-opt}$'$]
  Let $(L,b,\beta)$ be a gap-free SDP problem.
  Suppose there is a unique optimal solution $x$ which is generic over
  $\QQ[L,\beta]$
  in $S^n_{+r}$.
  Then there exists a strictly complementary pair $(x, \Omega)$.
\end{citingthm}

The proofs follow exactly the proofs of the versions given earlier.
For instance, in Theorem~\ref*{thm:sdp-opt}$'$, we apply
Theorem~\ref*{thm:urexp}$'$ to the cone $S^n_+$ (which is defined
over~$\QQ$) and the projection onto $S^n/(L + \langle \beta \rangle)$,
which is defined over whatever field is needed to define $L$ and
$\beta$.

\subsection{Relation to previous results}

A fundamental result 
of~\cite{alizadeh1997complementarity,pataki2001generic} 
on strict complementarity can be
summarized as follows.

\begin{otheorem}
\label{thm:AHOsc}
Suppose $(L,b,\beta)$ is a generic SDP program that is gap free.
Then  $(L,b,\beta)$ 
admits a strictly complementary pair of solutions.
\end{otheorem}

This result is neither stronger nor weaker than our Theorem~\ref{thm:sdp-opt}.
In particular  Theorem~\ref{thm:AHOsc} requires that all parameters
be generic.
In contrast, Theorem~\ref{thm:sdp-opt} does not assume genericity of any
of the parameters but rather assumes genericity of the solution within its
rank. 
Indeed, in our
application to rigidity of graphs, $\beta=0$, which is obviously not
generic.  (The parameter $b$ is also not usually generic.)
In fact,
there are very few situations where both theorems can
apply, as the following proposition shows.

\begin{proposition}
  In an SDP problem $(L,b,\beta)$ with $\dim L = D$ and primal
  solution~$x$ of rank~$r$, if
  Theorem~\ref{thm:AHOsc} applies, then
  \[
  \binom{n-r+1}{2} \le D.
  \]
  On the other hand, if Theorem~\ref{thm:sdp-opt} applies, then
  \[
  \binom{n-r+1}{2} \ge D.
  \]
\end{proposition}

\begin{proof}
  For a given $(L,b,\beta)$, if $b$ is generic
over  $\QQ[L]$,
 then $r$ satisfies the first inequality, as
  described 
in~\cite[Theorem 12]{alizadeh1997complementarity} 
and~\cite[Proposition 5]{nie2010algebraic}.
In particular, for generic $b$, the intersection of $L+b$ with 
$S^n_{+r}$ must be transversal, which from a dimension count gives
the first inequality.
  This takes care of the first part of the proposition.
  (If $\beta$ is generic over  $\QQ[L]$,
  as in Theorem~\ref{thm:AHOsc}, there is also
  an upper bound on~$r$:
  $\binom{r+1}{2} \leq \binom{n+1}{2}-D$.  But we do not need this.)

  For the second part, recall that if Theorem~\ref{thm:sdp-opt}
  applies, there is a 
point~$x$, generic
over
  $\QQ[L,\beta]$
  in $S^n_{+r}$
which is the unique solution
  to the SDP problem $(L, b, \beta)$.
  Recall that there is an associated feasibility problem
  $(L', b', 0)$.  Let $\pi: S^n \to S^n/L'$ be the associated
  projection and $\pi_+$ its restriction to $S^n_+$.
  Since $x$ is unique (in both $(L,b,\beta)$ and $(L',b',0)$), 
  $x$ is the only point in $S^n_+$ mapping to $\pi(x)$, and $\pi(x)
  \in \bdy \pi(S^n_+)$.
  Moreover, since $x$ is generic in $S^n_{+r}$, there is
  an open neighborhood $U$ of~$x$ in $S^n_{+r}$ with these
  properties.  In particular, $U$ maps injectively by $\pi$ to $\bdy
  \pi(S^n_+) \subset S^n/L'$.

  $S^n_{+r}$ has dimension
  $\binom{n+1}{2}-\binom{n-r+1}{2}$.  On the other hand, $\dim
  (S^n/L') \le\binom{n+1}{2}-D+1$ (with
  equality iff $\beta \not\in L^\perp$ or equivalently $L' \ne L$),
  so $\dim (\partial\pi(S^n_+))\le \binom{n+1}{2}-D$.  
  If there is a smooth injection from $S^n_{+r}$ to $\partial
  \pi(S^n_+)$, we must have $\binom{n-r+1}{2}\ge D$, as desired.
\end{proof}

\subsection{Complexity of universal rigidity}
\label{sec:comp}

Assuming $(\prho,\Gamma)$ is not at a conic at infinity and is
translated so that its affine span does not include~$0$,
$(\prho,\Gamma)$ is UR iff
there is no  higher rank solution to the SDP than 
$\prho$~\cite{Alfakih07:DimensionalRigidity}.
Thus, to test for universal rigidity,
the main step is to test if there is a feasible solution with rank
higher than that of a known input 
feasible
solution~$\prho$ (given say as integers).
Numerically speaking, SDP optimization algorithms that use 
interior point methods
produce approximate solutions of the highest possible 
rank~\cite{guler1993convergence,de1997initialization} 
and so in practice  one could try such a method to
produce a guess about the universal rigidity of $(\prho,\Gamma)$.

The complexity of getting a definitive answer is a trickier question,
even 
assuming 
one could reduce the UR question to one of ``yes-no'' SDP feasibility.
In particular, 
an approximate solution to an SDP optimization problem can be found
in polynomial time but 
the complexity of the  SDP feasibility problem remains
unknown, even with strict complementarity.
See~\cite{ramana1997exact} for formal details.

Theorem~\ref{thm:generic-ur-min-stress-kernel} tells us that for
generic inputs, the UR question can be answered
by finding the highest rank dual optimal solution.
(A framework with integer coordinates will not be generic;
however, if the integers are large enough we are likely to avoid all
special behavior, as in \cite[Section 5]{GHT10:GGR}.)
This does not appear to be any help in determining the complexity of 
testing algorithmically for 
universal rigidity.


\bibliographystyle{hamsalpha}
\bibliography{graphs,convex,alggeom}

\end{document}